\documentclass[11pt, oneside]{amsart}   	
\usepackage{geometry}                		
\usepackage{graphicx}				
\usepackage{comment, amssymb, amsmath, color}

\newtheorem{theorem}{Theorem}

\newtheorem{lemma}{Lemma}
\newtheorem{corollary}{Corollary}

\begin{document}

\title{Power partitions}
\author{Ayla Gafni}
\address{Pennsylvania State University \\ 109 McAllister Bldg \\ University Park, PA 16802 \\ gafni@math.psu.edu}

\begin{abstract}
In 1918, Hardy and Ramanujan published a seminal paper which included an asymptotic formula for the partition function.  In their paper, they also claim without proof an asymptotic equivalence for $p^k(n)$, the number of partitions of a number $n$ into $k$-th powers.  In this paper, we provide an asymptotic formula for $p^k(n)$, using the Hardy-Littlewood Circle Method.  We also provide a formula for the difference function $p^k(n+1)-p^k(n)$.   As a necessary step in the proof, we obtain a non-trivial bound on exponential sums of the form $\sum_{m=1}^q e(\frac{am^k}{q})$.
\end{abstract}

\maketitle

\section{Introduction}
A partition of number $n$ is a non-increasing sequence of positive integers whose sum is equal to $n$.  Fix an integer $k\ge2$.   Define $p^k(n)$ to be the number of partitions of $n$ in which all parts are perfect $k$-th powers.  This sequence has the generating function
$$\Psi_k(z) := \sum_{n=0}^\infty p^k(n) z^n = \prod_{n=1}^\infty (1 - z^{n^k})^{-1}.$$
In 1918, Hardy and Ramanujan \cite{Hardy1918} published a seminal paper introducing a new method for computing asymptotic formulae for integer sequences, called the Circle Method.  In their paper, they list a number of problems to which their methods can be applied.  In particular, they state (without proof) the following asymptotic equivalence for the number of partitions of $n$ into $k$-th powers: 
\begin{equation}\label{ramanujan}
\log p^k(n) \sim (k+1) \left(\frac{1}{k} \Gamma\left(1+\frac{1}{k}\right)\zeta\left(1+\frac{1}{k}\right)\right)^{k/(k+1)} n^{1/(k+1)}.
\end{equation}

 In 1934, E.~Maitland Wright \cite{Wright1934} gave a precise asymptotic formula for this restricted partition function. His proof requires a number of complicated objects including generalized Bessel functions.  In this paper, we provide a new asymptotic formula for the number of partitions into $k$-th powers, using a relatively simple implementation of the Hardy-Littlewood Circle Method.  The special case $k=2$ is treated by R.C. Vaughan \cite{Vaughan2015}.  This work is a generalization of Vaughan's result.

\begin{theorem}\label{main theorem} Let $n$ be a sufficiently large natural number, and choose positive numbers $X$ and $Y$ satisfying 
\begin{align*} n &= X\left(  \frac{1}{k^2}\, \zeta (\frac{k+1}{k})\Gamma(\frac{1}{k})X^{1/k} - \frac{1}{2}-\frac{1}{2}\zeta(-k)X^{-1}\right), \\ 
 Y  &= \frac{k+1}{2k^3} \,\zeta (\frac{k+1}{k})\Gamma(\frac{1}{k})X^{1/k} - \frac{1}{4}.
 \end{align*}
Then, for each $J\in\mathbb{N}$, there are real numbers $c_1, c_2, \ldots, c_{J}$ (independent of $n$), so that
\begin{equation}\label{main result}
p^k(n) = \frac{\exp\left(\frac{k+1}{k^2}\zeta(\frac{k+1}{k})\Gamma(\frac{1}{k})X^{\frac{1}{k}} - \frac{1}{2}\right)}{(2\pi)^{\frac{k+2}{2}}X^{\frac{3}{2}}Y^{\frac{1}{2}}} \left( \pi^{\frac{1}{2}} + \sum_{j=1}^{J} c_jY^{-j} + O(Y^{-J-1})\right).
\end{equation}
\end{theorem}

Note that 
$$X \sim \left(\frac{1}{k^2}\, \zeta (\frac{k+1}{k})\Gamma(\frac{1}{k})\right)^{-\frac{k}{k+1}} n^{\frac{k}{k+1}} , \quad Y \sim  \frac{k+1}{2k} \left(\frac{1}{k^2}\, \zeta (\frac{k+1}{k})\Gamma(\frac{1}{k})\right)^{\frac{k}{k+1}} n^{\frac{1}{k+1}}.$$
If we take the logarithm of both sides of (\ref{main result}), we obtain the asymptotic equivalence
$$\log p^k(n) \sim \frac{k+1}{k^2}\zeta(\frac{k+1}{k})\Gamma(\frac{1}{k})X^{\frac{1}{k}}.$$
Re-writing the right-hand side to be in terms of $n$ yields
$$\log p^k(n)\sim  (k+1)\left(\frac{1}{k^2}\, \zeta (\frac{k+1}{k})\Gamma(\frac{1}{k})\right)^{\frac{k}{k+1}} n^{\frac{1}{k+1}},$$
which is the same as Hardy and Ramanujan's claim (\ref{ramanujan}).

It is worthwhile to remark that it is possible to compute the coefficients $c_j$ that appear in (\ref{main result}).  However, their closed form is sufficiently complicated to make the statement of the theorem unreadable.  Thus we have omitted the closed form for $c_j$ from this paper.   Section \ref{computing the coefficients} provides an outline of how one might go about computing the values of $c_j$, given a particular choice of $J$.

The methods used to find the asymptotic formula in Theorem \ref{main theorem} can also be used to estimate the growth of $p^k(n)$.  This yields the following:

\begin{theorem}\label{difference theorem}  Let $n, X,$ and $Y$ be as above.  Then there are real numbers $d_1, d_2, \ldots, d_{k}$ (independent of $n$), so that
\begin{equation}\label{difference result}
p^k(n+1) - p^k(n) =\frac{\exp\left(\frac{k+1}{k^2}\zeta(\frac{k+1}{k})\Gamma(\frac{1}{k})X^{\frac{1}{k}} - \frac{1}{2}\right)}{(2\pi)^{\frac{k+2}{2}}X^{\frac{5}{2}}Y^{\frac{1}{2}}} \left( \pi^{\frac{1}{2}} + \sum_{j=1}^{k-1} d_jY^{-j} + O(Y^{-k})\right).
\end{equation}
\end{theorem}

From Theorem \ref{difference theorem} we can immediately deduce an asymptotic equivalence:

\begin{corollary} \label{difference cor} Let $n$ and $X$ be as above.  Then 
\begin{equation}\label{difference equiv}
p^k(n+1) - p^k(n) \sim \frac{p^k(n)}{X}
\end{equation}
as $n\rightarrow\infty$.
\end{corollary}

Before proceeding with the proof of the main result, we will need a few definitions.  Let
$$\Phi_k(z) := \sum_{j=1}^\infty\sum_{n=1}^\infty \frac{1}{j} z^{jn^k}.$$
We then have 
$$\Psi_k(z) = \exp(\Phi_k(z)).$$
For convenience, we also define
\begin{equation}\label{Delta}
\rho := e^{-1/X} \quad \mbox{and}\quad  \Delta := (1+4\pi^2 X^2 \Theta^2)^{-1/2}
\end{equation}
for  $X, \Theta\in \mathbb{R}$ with $X\ge1$.

We will prove Theorem \ref{main theorem} using the Hardy-Littlewood Method.  In dealing with the major arcs, we require a bound on exponential sums of the form
\begin{equation}
S_k(r,b) := \sum_{m=1}^r e\left(\frac{b m^k}{r}\right).
\end{equation}
Here we use the standard notation $e(\alpha) = e^{2\pi i \alpha}$.  The proof of the bound is quite simple, yet doesn't appear to exist in the literature.  We state it here as a lemma, as it may be of interest beyond the scope of this paper.
\newline

\begin{lemma} \label{exponential sum lemma} For each $k\ge2$ there exists a positive constant $\delta_k$ such that $|S_k(r,b)| \le (1-\delta_k)r$ for all $r>1$ and $b\in\mathbb{Z}$ with $(r,b)=1.$
\end{lemma}
{\bf Proof}~  By Theorem 4.2 of \cite{Vaughan1997} we have that $S_k(r,b) \ll r^{1-1/k}$ for $(r,b)=1$.  Thus there exists $C_k$ such that $|S_k(r,b)| \le C_k r^{1-1/k}$.  So, we can find $R$ sufficiently large and $\nu_k>0$ such that $|S_k(r,b)| \le (1-\nu_k) r$ for all  $r\ge R, b\in\mathbb{Z}$ with $(r,b)=1$.

 If $1<r<R$, then there is at least one term in $S_k(r,b)$ that is not equal to $1$.  This term is of the form $e(bm^k/r)$.  Therefore 
 $$S_k(r,b) \le |r - 1 + e(bm^k/r)| \le |r - 1 + e(1/r)| \le |r - 1 + e(1/R)| < (1-\eta_k)r,$$
 for some $\eta_k>0$.
 
Let $\delta_k = \min(\nu_k, \eta_k).$  Then $|S_k(r,b)| \le (1-\delta_k)r$ for all $r>1,\  b\in\mathbb{Z}$ with $(r,b)=1.$ \hfill $\square$

\section{Auxiliary Lemmas}

At several points in the proof of Theorem \ref{main theorem}, we will need to estimate the value of $\Phi_k(\rho e(\Theta))$.  So, before proving the theorem, we introduce two estimates for this expression.  Lemma \ref{near origin lemma} provides a very precise estimate that will be used for $|\Theta|\le\frac{3}{8\pi X}$ and will establish the main term of (\ref{main result}).  Lemma \ref{major arcs lemma} provides a less precise estimate that will be used to deal with the major arcs (excluding $|\Theta|\le\frac{3}{8\pi X}$ ).   The estimate needed for the minor arcs is provided in the body of the proof of Theorem \ref{main theorem}.
\newline

\begin{lemma} \label{near origin lemma} Suppose $\Theta\in\mathbb{R}$ and $X\ge 1$. 
If $X\Delta^3 \ge 1$, then 
\begin{align}\Phi_k(\rho e(\Theta)) & = \frac{1}{k}\, \zeta \left(\frac{k+1}{k}\right)\Gamma\left(\frac{1}{k}\right)\left(\frac{X}{1-2\pi i X \Theta}\right)^{1/k} - \frac{1}{2}\log\left(\frac{(2\pi)^k X}{1-2\pi i X \Theta}\right) 
\\ &+\frac{1}{2}\zeta(-k)\left(\frac{1-2\pi i X \Theta}{X}\right) 
+ O\left(\Delta^{-1/2}\exp\left(-\frac{1}{k}\left(2(\pi \Delta)^{k+1}X\right)^{1/k} \right)\right).\nonumber  \end{align}
\end{lemma}

\noindent{\bf Proof}~  We have
$$\Phi_k(\rho e(\Theta)) = \sum_{j=1}^\infty\sum_{n=1}^\infty \frac{1}{j} \exp\left(-jn^k\left(\frac{1}{X} - 2\pi i \Theta\right)\right).$$ 
Using a Mellin Transform (Theorem C.4 of \cite{Montgomery2006}), this is equal to 
$$\frac{1}{2\pi i}\sum_{j=1}^\infty\sum_{n=1}^\infty \frac{1}{j} \int_{c-i\infty}^{c+i\infty}\Gamma(s)j^{-s} n^{-ks}\left(\frac{X}{1 - 2\pi i X\Theta}\right)^s\,ds$$
for $c>0$.  Here $\left(\frac{X}{1-2\pi i X \Theta}\right)^{s} $ denotes $\exp\left(s \log\left(\frac{X}{1-2\pi i X \Theta}\right)\right)$ where the logarithm is defined by continuous variation of $ \log\left(\frac{X}{1-2\pi i X \theta}\right)$ as $\theta$ varies continuously from $0$ to $\Theta$  through real values.  The series $\zeta(s+1)$ and $\zeta(ks)$ converge absolutely and uniformly for $\mathfrak{R}s\ge 1/k + \delta$ where $\delta$ is any positive number.  Hence, for any real $c > 1/k$ we have,
$$\Phi_k(\rho e(\Theta)) = \frac{1}{2\pi i} \int_{c-i\infty}^{c+i\infty} \zeta(s+1)\zeta(ks)\left(\frac{X}{1-2\pi i X \Theta}\right)^{s} \Gamma(s) \,ds.$$
 
Since $\Gamma(\frac{1}{k})$ is well-defined, the integrand has a simple pole at $s = 1/k$ with residue
$$\frac{1}{k}\, \zeta (\frac{k+1}{k})\Gamma(\frac{1}{k})\left(\frac{X}{1-2\pi i X \Theta}\right)^{1/k}.$$
The integrand also has a double pole at $s=0$ from $\zeta(s+1)\Gamma(s)$.  The Laurent expansion of $\zeta(s+1)\Gamma(s)$ at $s=0$ is of the form 
$$\frac{1}{s^2} + \sum_{j=0}^{\infty} a_j s^j,$$
so the residue of the integrand at $s=0$ is 
$$\zeta(0) \log\left(\frac{X}{1-2\pi i X \Theta}\right) + k\zeta'(0).$$
We recall (see, for example, \cite{Montgomery2006}) that $\zeta(0) = -\frac{1}{2}$ and $\zeta'(0) = -\frac{1}{2}\log 2\pi$, so the residue of the integrand at $s=0$ is
$$- \frac{1}{2}\log\left(\frac{(2\pi)^k X}{1-2\pi i X \Theta}\right).$$

The $\Gamma$-function has simple poles at the negative integers, but the $\zeta$-function also has zeros at the negative even integers.  If $s\le -2$ is an integer, then either $ks$ or $s+1$ is an even integer, so the poles of $\Gamma$ are cancelled by the zeros of $\zeta(ks)\zeta(s+1)$.  This leaves one simple pole at $s=-1$.  The residue here is 
$$-\zeta(0)\zeta(-k)\left(\frac{1-2\pi i X \Theta}{X}\right)  = \frac{1}{2}\zeta(-k)\left(\frac{1-2\pi i X \Theta}{X}\right).$$

By the functional equation for the zeta function in its asymmetrical form, we have
\begin{align*}& \zeta(s+1) \zeta(ks) \Gamma(s) 
\\ &= 2^{(k+1)s+1} \pi^{(k+1)s-1}  \cos(\frac{\pi s}{2})\sin(\frac{k\pi s}{2})\zeta(-s) \zeta(1-ks) \Gamma(-s) \Gamma(s)\Gamma(1-ks).
\end{align*}
Moreover, by the reflection formula for the gamma function, 
$$\Gamma(-s) \Gamma(s) = \frac{\pi}{-s \sin(\pi s)}.$$
Combining this with the fact that $\Gamma(1- ks) = -ks\Gamma(-ks)$, we have
$$\zeta(s+1) \zeta(ks) \Gamma(s) = (2k)(2\pi)^{(k+1)s} \zeta(-s) \zeta(1-ks) \Gamma(-ks) \frac{\cos(\frac{\pi s}{2})\sin(\frac{k\pi s}{2})}{\sin(\pi s)}.$$
Note that 
$$\frac{\cos(\frac{\pi s}{2})\sin(\frac{k\pi s}{2})}{\sin(\pi s)} = \frac{\sin(k\pi s/2)}{2\sin(\pi s/2)} \ll e^{(k-1)|t|\pi/2}.$$
Therefore, by Stirling's formula, when $\sigma = \mathfrak{R}s\le -3/2$, 
$$\zeta(s+1) \zeta(ks) \Gamma(s) \ll (2\pi)^{(k+1)\sigma} k^{-k\sigma} |s|^{-1/2-k\sigma} e^{-\pi|t|/2}.$$

From the argument in \cite{Vaughan2015}, we see that
$$\left|\left(\frac{X}{1-2\pi i X \Theta}\right)^{s}\right| \le (X\Delta)^\sigma \exp\left(|t|\left(\frac{\pi}{2}-\Delta\right)\right).$$
Therefore the integrand is 
$$\ll (2\pi)^{(k+1)\sigma} k^{-k\sigma} |s|^{-1/2-k\sigma}  (X\Delta)^\sigma e^{-\Delta|t|}.$$

Let $R\ge\frac{3}{2}$, and move the vertical line of integration to the line $\mathfrak{R}s = -R$.  On this line the integrand is 
$$\ll \left(\frac{1}{(2\pi)^{(k+1)}X\Delta}\right)^R  k^{kR} |R+it|^{-1/2+kR}  e^{-\Delta|t|}.$$
On the pieces with $|t|>R$, the integrand is 
$$\ll \left(\frac{1}{2 \pi^{(k+1)}X\Delta}\right)^R  k^{kR} |t|^{-1/2+kR}  e^{-\Delta|t|}.$$
After a change of variable ($y=\Delta |t|$), this contributes
$$\ll \left(\frac{k^k}{2 (\pi\Delta)^{(k+1)}X}\right)^R  \Delta^{-1/2} \Gamma(kR+1/2),$$
which by Stirling's formula is 
$$\ll \left(\frac{k^{2k}e^{-k}}{2 (\pi\Delta)^{(k+1)}X}\right)^R  R^{kR} \Delta^{-1/2} .$$
Meanwhile, the part of the integral with $|t|\le R$ contributes
$$\ll \left(\frac{1}{2\pi^{(k+1)}X\Delta}\right)^R k^{kR} R^{-1/2+kR} \Delta^{-1}.$$
Combining the estimates we see that the integral is 
$$\ll \left(\frac{k^{2k}}{2 e^{k}(\pi\Delta)^{(k+1)}X}\right)^R  R^{kR} \Delta^{-1/2} .$$

The choice of $R$ which minimizes the expression above is 
$$R_0 = \left(\frac{2 (\pi\Delta)^{(k+1)}X}{k^{2k}}\right)^{1/k}.$$
Let $R = \max(R_0,3/2)$.  Then the integral is 
$$\ll \exp\left(-\frac{1}{k}\left(2 (\pi\Delta)^{(k+1)}X\right)^{1/k}\right) \Delta^{-1/2}.$$

Applying the Cauchy residue theorem, we obtain the result.
\hfill $\square$
\newline
\newline

\begin{lemma} \label{major arcs lemma} Suppose that $X\in\mathbb{R}, X>1, \Theta\in\mathbb{R}, a\in\mathbb{Z}, q\in\mathbb{N}, (a,q) = 1$ and $\theta = \Theta-a/q$.  Then
\begin{align*}
\Phi_k(\rho e(\Theta)) & = \Gamma\left(\frac{k+1}{k}\right)\left(\frac{X}{1-2\pi i X\Theta}\right)^{1/k}\sum_{j=1}^\infty \frac{S_k(q_j, a_j)}{j^{\frac{k+1}{k}}q_j} \\&
+ O(q^{1/2+\varepsilon}\log{X}(1+X^{1/2}|\theta|^{1/2}))
\end{align*}
where $q_j = q/(q,j)$ and $a_j = aj/(q,j)$.
\end{lemma}

\noindent{\bf Proof}~  From the definition, we have
$$\Phi_k(\rho e(\Theta)) =  \sum_{j=1}^\infty\sum_{n=1}^\infty \frac{1}{j} e^{-n^kj/X}e(jn^k\Theta).$$
We write
$$e^{-n^kj/X} = \int_n^\infty k x^{k-1} jX^{-1} e^{-x^kj/X} \,dx.$$
Thus 
$$\Phi_k(\rho e(\Theta)) =  \sum_{j=1}^\infty \frac{1}{j}  \int_0^\infty k x^{k-1} jX^{-1} e^{-x^kj/X} \sum_{n\le x} e(jn^k\Theta)\,dx.$$
It is useful to observe the crude bound 
$$\int_0^\infty k x^{k-1} jX^{-1} e^{-x^kj/X} \sum_{n\le x} e(jn^k\Theta)\,dx \ll \int_0^\infty k x^{k} jX^{-1} e^{-x^kj/X}\,dx.$$
Using integration by parts, this is 
$$\int_0^\infty e^{-x^kj/X}\,dx = \left(\frac{X}{j}\right)^{1/k}\int_0^\infty e^{-y^k}\,dx \ll  \left(\frac{X}{j}\right)^{1/k}.$$ 
Let $J$ be a parameter at our disposal.  Then 
$$ \sum_{j=J+1}^\infty \frac{1}{j}  \int_0^\infty k x^{k-1} jX^{-1} e^{-x^kj/X} \sum_{n\le x} e(jn^k\Theta)\,dx \ll \sum_{j=J+1}^\infty \frac{1}{j}\left(\frac{X}{j}\right)^{1/k} \ll \left(\frac{X}{J}\right)^{1/k}.$$

It remains to consider 
\begin{equation} \sum_{j=1}^J \frac{1}{j}  \int_0^\infty k x^{k-1} jX^{-1} e^{-x^kj/X} \sum_{n\le x} e(jn^k\Theta)\,dx.\label{finite sum}\end{equation}
By Theorem 4.1 of \cite{Vaughan1997}, we have
$$ \sum_{n\le x} e(jn^k\Theta) = q_j^{-1} S_k(q_j, a_j)\int_0^x e(j\theta\gamma^k)\,d\gamma + O\left(q_j^{1/2+\varepsilon}(1+x^k j|\theta|)^{1/2}\right).$$
Hence the expression in (\ref{finite sum}) is equal to 
\begin{equation} \sum_{j=1}^J \frac{S_k(q_j, a_j) }{jq_j} \int_0^\infty k x^{k-1} jX^{-1} e^{-x^kj/X} \int_0^x e(j\theta\gamma^k)\,d\gamma\,dx + E_1,\label{after thm 4.1}\end{equation}
where
$$E_1 \ll \sum_{j=1}^J \frac{q_j^{1/2+\varepsilon}}{j}  \int_0^\infty k x^{k-1} jX^{-1} e^{-x^kj/X}(1+(x^k j|\theta|)^{1/2}) \,dx.$$
By integration by parts we have
$$E_1 \ll \sum_{j=1}^J \frac{q_j^{1/2+\varepsilon}}{j} \left(1+ \frac{k}{2}(j|\theta|)^{1/2}\int_0^\infty x^{k/2-1} e^{-x^k j/X}\,dx\right).$$
We now make the substitution $y = x^k jX^{-1}$. Then we have $x =( yX/j)^{1/k}$, and $dx = \frac{1}{k}y^{1/k-1}(X/j)^{1/k}dy$.  So the integrand becomes 
$$\left(\frac{yX}{j}\right)^{1/2-1/k}e^{-y}\ \frac{1}{k}y^{1/k-1}\left(\frac{X}{j}\right)^{1/k} = \frac{1}{k}\left(\frac{X}{j}\right)^{1/2} y^{-1/2} e^{-y}.$$
Thus
\begin{align*}
E_1 & \ll \sum_{j=1}^J \frac{q_j^{1/2+\varepsilon}}{j} \left(1+ \frac{1}{2}(X|\theta|)^{1/2}\int_0^\infty y^{-1/2} e^{-y}\,dy\right) \\ &
 \ll q^{1/2+\varepsilon}\log{J}\left(1+ |\theta|^{1/2}X^{1/2}\right).
\end{align*}

We now turn our attention to the main term of the expression in (\ref{after thm 4.1}).  By integration by parts, this is 
$$ \sum_{j=1}^J \frac{S_k(q_j, a_j) }{jq_j} \int_0^\infty e^{-x^k j/X}e(x^k j\theta)\,dx.$$
The integral here is 
$$\int_0^\infty exp(-x^k jX^{-1}(1-2\pi i X\theta))\,dx.$$
We make the substitution $z = \left(jX^{-1}(1-2\pi i X\theta)\right)^{1/k} x$.  Choose $\phi$ so that $|\phi|<\pi$ and 
$$\frac{1-2\pi i X\theta}{|1-2\pi i X\theta|} = e^{i\phi}.$$
We thus obtain 
$$z = \left(jX^{-1}|1-2\pi i X\theta|\right)^{1/k} e^{i\phi/k} x.$$
This gives
$$\int_0^\infty exp(-x^k jX^{-1}(1-2\pi i X\theta))\,dx =  \left(\frac{X}{j(1-2\pi i X\theta)}\right)^{1/k}\int_\mathcal{L} e^{-z^k}\,dz,$$
where $\mathcal{L}$ is the ray $\{z = x e^{i\phi/k}:  0\le x < \infty\}$.   By Cauchy's theorem, the integral here is $\Gamma(\frac{k+1}{k})$.  Putting everything together, we have
\begin{align*}
\Phi_k(\rho e(\Theta)) & = \Gamma\left(\frac{k+1}{k}\right)\left(\frac{X}{1-2\pi i X\theta}\right)^{1/k}\sum_{j=1}^J \frac{S_k(q_j, a_j) }{j^{\frac{k+1}{k}}q_j}
\\ & + O\left(q^{1/2+\varepsilon}\log{J}\left(1+ |\theta|^{1/2}X^{1/2}\right) + \left(\frac{X}{J}\right)^{1/k}\right).
\end{align*}
Since $|S_k(q_j, a_j)| \le q_j$, we have 
$$ \Gamma\left(\frac{k+1}{k}\right)\left(\frac{X}{1-2\pi i X\theta}\right)^{1/k}\sum_{j=J+1}^\infty \frac{S_k(q_j, a_j) }{j^{\frac{k+1}{k}}q_j} \ll \left(\frac{X}{J}\right)^{1/k}.$$
Thus we can extend the sum in the main term to infinity.  Setting $J=X$ gives
\begin{align*}
\Phi_k(\rho e(\Theta))&  = \Gamma\left(\frac{k+1}{k}\right)\left(\frac{X}{1-2\pi i X\theta}\right)^{1/k}\sum_{j=1}^\infty \frac{S_k(q_j, a_j) }{j^{\frac{k+1}{k}}q_j} \\ &
 + O\left(q^{1/2+\varepsilon}\log{X}\left(1+ |\theta|^{1/2}X^{1/2}\right) \right),
 \end{align*}
as desired.
\hfill $\square$
\newline
\newline

\section{Proof of Theorem \ref{main theorem}}

We prove the theorem using the Hardy-Littlewood circle method.  From Cauchy's theorem, we have
$$p^k(n) = \int_{0}^{1} \rho^{-n} \exp(\Phi_k (\rho e(\Theta))-2\pi i n \Theta)\, d\Theta.$$
By the periodicity of the integrand, we may replace the unit interval by any interval of length 1.  We will use the interval $\mathcal{U} = (-X^{\frac{1}{k}-1}, 1- X^{\frac{1}{k}-1}]$.  This is a convenient choice because the main contribution to the integral comes from $\Theta$ near the origin.  Using $\mathcal{U}$ instead of $[0,1]$ prevents that region from being split in two.

For $a,q\in\mathbb{N}$ with $(a,q)=1$, define
$$\mathfrak{M}(q,a) = \{ \Theta\in \mathcal{U} : |\Theta - \frac{a}{q}|\le q^{-1}X^{\frac{1}{k}-1}\},$$
and let 
$$\mathfrak{M} = \bigcup_{1\le a\le q \le X^{1/k}} \mathfrak{M}(q,a).$$
We refer to these disjoint intervals as the major arcs, and we define the minor arcs to be the complement of the major arcs, namely
$$\mathfrak{m} =\mathcal{U} \setminus\mathfrak{M}.$$

In a typical implementation of the circle method, one would split the intergral into the major arcs and the minor arcs, with the major arcs making up the main term of the asymptotic formula.  However, in our case, the contribution from $\mathfrak{M}(1,0)$ is significantly greater than the contribution from the rest of the major arcs.  So we will split the integral into three main parts, namely
$$p^k(n) = \left\{\int_{\mathfrak{M}(1,0)} + \int_{\mathfrak{M}\setminus\mathfrak{M}(1,0)} + \int_{\mathfrak{m}} \right\} \rho^{-n} \exp(\Phi_k (\rho e(\Theta))-2\pi i n \Theta)\, d\Theta.$$
We will treat $\mathfrak{M}\setminus\mathfrak{M}(1,0)$ and $\mathfrak{m}$ in the way one would traditionally treat the major and minor arcs, respectively, but the major arcs with $q>1$ will not contribute to the main term of the asymptotic formula.  Rather they will be ``thrown away'' into the error term.  The main term of the asymptotic formula will come from the first part of the integral, when $\Theta$ is close to the origin.  We examine that piece first.  

We first consider 
\begin{equation}\label{k int near 0} \int_{-3/(8\pi X)}^{3/(8\pi X)} \rho^{-n} \exp(\Phi_k (\rho e(\Theta))-2\pi i n \Theta)\, d\Theta.
\end{equation}
When $|\Theta|\le1/X$,
$$(1+4\pi^2)^{-\frac{1}{2}} \le \Delta \le 1$$ 
and by Lemma \ref{near origin lemma}  we have 
\begin{align}\Phi_k(\rho e(\Theta)) & = \frac{1}{k}\, \zeta \left(\frac{k+1}{k}\right)\Gamma\left(\frac{1}{k}\right)\left(\frac{X}{1-2\pi i X \Theta}\right)^{1/k} \\ &\nonumber
- \frac{1}{2}\log\left(\frac{(2\pi)^k X}{1-2\pi i X \Theta}\right) 
+\frac{1}{2}\zeta(-k)\left(\frac{1-2\pi i X \Theta}{X}\right) \\ &
+ O\left(\Delta^{-1/2}\exp\left(-\frac{1}{k}\left(2(\pi \Delta)^{k+1}X\right)^{1/k} \right)\right).\nonumber  \end{align}
Thus we have 
\begin{align}\label{k phi to xi}\rho^{-n} &\Psi_k(\rho e(\Theta)) \\ &
= \rho^{-n} \exp(\Xi_k(\rho e(\Theta)))\left(1+O\left(\Delta^{-1/2}\exp\left(-\frac{1}{k}\left(2(\pi \Delta)^{k+1}X\right)^{1/k} \right)\right)\right)
\nonumber
\end{align}
where
\begin{align*}\Xi_k(\rho e(\Theta)) & = \frac{1}{k}\, \zeta \left(\frac{k+1}{k}\right)\Gamma\left(\frac{1}{k}\right)\left(\frac{X}{1-2\pi i X \Theta}\right)^{1/k} \\ &
- \frac{1}{2}\log\left(\frac{(2\pi)^k X}{1-2\pi i X \Theta}\right) 
+\frac{1}{2}\zeta(-k)\left(\frac{1-2\pi i X \Theta}{X}\right).
\end{align*}

We also write 
$$\frac{ X}{1-2\pi i X \Theta} = X \Delta e^{i\phi}$$
where $\phi =  \arg(1+2\pi i X \Theta)$.    Note that $0<|\phi|\le\pi/2$, so $0<\cos(\phi/k)<1$.  Hence
$$\left| \left( \frac{ X}{1-2\pi i X \Theta}\right)^{\frac{1}{k}}\right| = (X \Delta)^{\frac{1}{k}}.$$  So the  $O$-term in (\ref{k phi to xi}) contributes
\begin{equation}\label{ftn of Delta}\ll X^{-1/2} \exp\left( \frac{n}{X} + \frac{1}{k} \zeta (\frac{k+1}{k})\Gamma(\frac{1}{k})(X\Delta)^{\frac{1}{k}} 
+ \frac{1}{2} \zeta(-k)(X\Delta)^{-1}   - \frac{1}{k}\left(2(\pi\Delta)^{k+1} X\right)^{\frac{1}{k}}\right).\end{equation}
As a function of $\Delta$, the expression here has a unique local maximum at 
$$\Delta_0 =\frac{ \zeta (\frac{k+1}{k})\Gamma(\frac{1}{k})}{(k+1)2^{\frac{1}{k}}\pi^{\frac{k+1}{k}}}.$$
Direct computation shows that $\Delta_0 <1$ for $k\le 5$ and $\Delta_0 >1$ for $k\ge6$.   We will address these cases separately.  It will help to recall that 
$$\frac{n}{X} = \frac{1}{k^2} \zeta(\frac{k+1}{k})\Gamma(\frac{1}{k})X^{1/k}-\frac{1}{2} - \frac{1}{2}\zeta(-k)X^{-1}.$$

If $k\ge6$, then $\Delta_0>1$, so the expression in (\ref{ftn of Delta}) is monotonically increasing in $\Delta$ for $\Delta \le1$.  Thus we may replace $\Delta$ by 1 to obtain that  (\ref{ftn of Delta}) is
$$\ll  X^{-1/2} \exp\left(\frac{k+1}{k^2}\zeta(\frac{k+1}{k})\Gamma(\frac{1}{k})X^{\frac{1}{k}} - \frac{1}{2} - \frac{1}{k}\left(2\pi^{k+1} X\right)^{\frac{1}{k}}\right).$$
If $k=2$ or $k=4$, then $\zeta(-k)=0$ and $0<\Delta_0<1$, so the expression in (\ref{ftn of Delta}) is 
$$ \ll  X^{-1/2} \exp\left(\frac{k+1}{k^2}\zeta(\frac{k+1}{k})\Gamma(\frac{1}{k})X^{\frac{1}{k}} - \frac{1}{2} - \frac{1}{k}\left(2(\pi\Delta_0)^{k+1} X\right)^{\frac{1}{k}}\right).$$
If $k=3$ or $k=5$, then $\frac{2}{5}<\Delta_0<1$.  Note that $\zeta(-3) = \frac{1}{120}, \zeta(-5) = \frac{-1}{252}$.  Thus the terms involving $\zeta(-k)$ in  (\ref{ftn of Delta}) have absolute value less than 1, and certainly less than $|\frac{1}{k}\left(2(\pi\Delta)^{k+1} X\right)^{\frac{1}{k}}|$.  So, there exists some fixed $\delta>0$ such that for any $k\ge2$, the expression in (\ref{ftn of Delta}) is 
$$\ll  X^{-1/2} \exp\left(\frac{k+1}{k^2}\zeta(\frac{k+1}{k})\Gamma(\frac{1}{k})X^{\frac{1}{k}} - \frac{1}{2} -\delta X^{\frac{1}{k}}\right).$$
Hence the integral in (\ref{k int near 0}) is 
\begin{align}\label{k int near 0 w error} 
& \int_{-3/(8\pi X)}^{3/(8\pi X)}\rho^{-n} \exp(\Xi_k(\rho e(\Theta))-2\pi i n \Theta)\, d\Theta 
\\ & + O\left(X^{-1/2} \exp\left(\frac{k+1}{k^2}\zeta(\frac{k+1}{k})\Gamma(\frac{1}{k})X^{\frac{1}{k}} - \frac{1}{2} -\delta X^{\frac{1}{k}}\right)\right). \nonumber
\end{align}

We turn our attention to the main term in (\ref{k int near 0 w error}).  When $|\Theta| < 1/(2\pi X)$, the function 
$$F(\Theta)  := \Xi_k(\rho e(\Theta))-2\pi i n \Theta$$
can be expanded as
$$F(\Theta)  = \frac{1}{k} \zeta (\frac{k+1}{k})\Gamma(\frac{1}{k})X^{\frac{1}{k}}+\frac{1}{2}\zeta(-k)X^{-1} - \frac{1}{2}\log\left((2\pi)^k X\right) - Y(2\pi X \Theta)^2 + G(\Theta) $$
where
\begin{equation}\label{G defn}
G(\Theta) =  \sum_{j=3}^\infty \left( \frac{\zeta (\frac{k+1}{k})\Gamma(j+\frac{1}{k})X^{\frac{1}{k}} }{j!\ k} - \frac{1}{2j}\right) (2\pi i X \Theta)^j.
\end{equation}
Thus the integral in the main term of (\ref{k int near 0 w error}) becomes
\begin{align}
&\frac{\exp\left(\frac{n}{X} + \frac{1}{k} \zeta (\frac{k+1}{k})\Gamma(\frac{1}{k})X^{\frac{1}{k}}+\frac{1}{2}\zeta(-k)X^{-1} \right)}{(2\pi)^{\frac{k}{2}}X^{\frac{1}{2}}} \int_{-3/(8\pi X)}^{3/(8\pi X)} \exp\left(-Y(2\pi X \Theta)^2 + G(\Theta) \right)\, d\Theta \nonumber
\\& = \frac{\exp\left(\frac{k+1}{k^2}\zeta(\frac{k+1}{k})\Gamma(\frac{1}{k})X^{\frac{1}{k}} - \frac{1}{2}\right)}{(2\pi)^{\frac{k}{2}}X^{\frac{1}{2}}}\int_{-3/(8\pi X)}^{3/(8\pi X)} \exp\left(-Y(2\pi X \Theta)^2 + G(\Theta) \right)\, d\Theta. \label{k int near 0 w G}
\end{align}

We rewrite the coefficients in the series for $G(\Theta)$ as
$$ \frac{\zeta (\frac{k+1}{k})\Gamma(j+\frac{1}{k})X^{\frac{1}{k}} }{j!\ k} - \frac{1}{2j} = a_j Y + b_j$$
where 
$$ a_j = \left(\frac{2k^2}{k+1} \right)\frac{\Gamma(j+\frac{1}{k}) }{j!\ \Gamma(\frac{1}{k})}, \quad b_j = \frac{1}{2}\left(\frac{2k^2}{k+1} \right)\frac{\Gamma(j+\frac{1}{k}) }{j!\ \Gamma(\frac{1}{k})} - \frac{1}{2j}.$$  
We then rewrite the integral on the right-hand side of (\ref{k int near 0 w G}) as 
\begin{align*}
\int_{0}^{3/(8\pi X)}\left( \exp( G(\Theta)) +  \exp(G(-\Theta))\right)\exp\left(-Y(2\pi X \Theta)^2 \right)\, d\Theta \\ = \mathfrak{R} \int_{0}^{3/(8\pi X)} 2\exp\left(G(\Theta) -Y(2\pi X \Theta)^2 \right) \, d\Theta.
\end{align*}
We now make the change of variable $\phi = (2\pi X\Theta)^2 Y$, so that the right-hand side of (\ref{k int near 0 w G}) becomes
\begin{equation}\label{k int near 0 with H}
\frac{\exp\left(\frac{k+1}{k^2}\zeta(\frac{k+1}{k})\Gamma(\frac{1}{k})X^{\frac{1}{k}} - \frac{1}{2}\right)}{(2\pi)^{\frac{k+2}{2}}X^{\frac{3}{2}}Y^{\frac{1}{2}}} \int_{0}^{9Y/16} \mathfrak{R} \left(\exp\left(-\phi+H(\phi) \right)\right)\phi^{-1/2} \, d\phi
\end{equation}
where 
\begin{equation}\label{H defn}
H(\phi) = \sum_{j=3}^\infty i^j (a_j + b_j Y^{-1}) \phi^{\frac{j}{2}} Y^{1-\frac{j}{2}}.
\end{equation}

We have
 $$\mathfrak{R} \left(\exp\left(-\phi+H(\phi)\right) \right) \le e^{-\phi}\left|e^{H(\phi)}\right| = e^{-\phi}\exp\left(\mathfrak{R}H(\phi)\right)$$
and
$$\mathfrak{R}H(\phi) = Y \sum_{j=2}^\infty  (-1)^j (a_{2j} + b_{2j} Y^{-1}) \phi^{j} Y^{-j}.$$
For $j\ge2$, we have $0<a_{2j}\le a_4 = \frac{6k^2+5k+1}{12k^2}$ and $0\le b_{2j} \le a_{2j}$.  Hence, when $0\le\phi\le9Y/16$, we have
\begin{align*}
\left|\mathfrak{R}H(\phi) \right| & \le\left(\frac{6k^2+5k+1}{12k^2}\right) Y(1+ Y^{-1}) \sum_{j=2}^\infty \left( \frac{\phi}{Y}\right)^{j} 
\\& = \left(\frac{6k^2+5k+1}{12k^2}\right) \frac{Y(1+ Y^{-1})(\phi/Y)^2}{1-\frac{\phi}{Y}} 
\\& \le \frac{9}{7}\left(\frac{6k^2+5k+1}{12k^2}\right) (1+ Y^{-1})\phi <\left(1- \frac{1}{56k^2}\right)\phi,
\end{align*}
since $n$ is sufficiently large.  Hence for $Z>0$,
\begin{align*}\int _Z^{9Y/16} \mathfrak{R} \left(\exp\left(-\phi+H(\phi) \right)\right)\phi^{-1/2} \, d\phi & \le \int _Z^{9Y/16} e^{-\phi} \exp\left(\left(1- \frac{1}{56k^2}\right) \phi \right) \phi^{-1/2} \, d\phi
\\ & \ll Z^{-1/2}\int_Z^\infty e^{-\phi/(56k^2)}\,d\phi \ll  Z^{-1/2}e^{-Z/(56k^2)}.
\end{align*}
Take $Z = 56k^2 J \log Y$.  Then 
$$\int _Z^{9Y/16} \mathfrak{R} \left(\exp\left(-\phi+H(\phi) \right)\right)\phi^{-1/2} \, d\phi \ll Y^{-J}.$$
The integral in (\ref{k int near 0 w error}) now becomes
\begin{equation}\label{k int near 0 w Z}
\frac{\exp\left(\frac{k+1}{k^2}\zeta(\frac{k+1}{k})\Gamma(\frac{1}{k})X^{\frac{1}{k}} - \frac{1}{2}\right)}{(2\pi)^{\frac{k+2}{2}}X^{\frac{3}{2}}Y^{\frac{1}{2}}} \left( \int _0^Z \mathfrak{R} \left(\exp\left(-\phi+H(\phi) \right)\right)\phi^{-1/2} \, d\phi + O(Y^{-J}) \right).
\end{equation}

When $0\le\phi\le Z$, we have
$$\sum_{j=2J+3}^\infty  i^j (a_j + b_j Y^{-1}) \phi^{\frac{j}{2}} Y^{1-\frac{j}{2}} \ll \frac{\phi^{J+\frac{3}{2}}Y^{-\frac{1}{2}-J}}{1-(\phi/Y)^\frac{1}{2}} \ll \phi^{J+\frac{3}{2}}Y^{-\frac{1}{2}-J}.$$
Hence,
$$\exp\left(\sum_{j=2J+3}^\infty  i^j (a_j + b_j Y^{-1}) \phi^{\frac{j}{2}} Y^{1-\frac{j}{2}} \right) = 1+ O(\phi^{J+\frac{3}{2}}Y^{-\frac{1}{2}-J}).$$
Let 
$$H_J(\phi) = \sum_{j=3}^{2J+2}  i^j (a_j + b_j Y^{-1}) \phi^{\frac{j}{2}} Y^{1-\frac{j}{2}}.$$
Then
\begin{align*}
&\int _0^Z \mathfrak{R} \left(\exp\left(-\phi+H(\phi) \right)\right)\phi^{-1/2} \, d\phi
\\& = \int _0^Z \mathfrak{R} \left(\exp\left(-\phi+H_J(\phi) \right)\right)\left( 1+ O(\phi^{J+\frac{3}{2}}Y^{-\frac{1}{2}-J})\right) \phi^{-1/2}\, d\phi.
\end{align*}
Similar to the argument for $H$, we have $\mathfrak{R}H_J(\phi) < \left(1-\frac{1}{56k^2}\right)\phi$ and so the error term here contributes
$$\ll  Y^{-\frac{1}{2} -J}\int_0^\infty e^{-\phi/(56k^2)} \phi^{J+1}\,d\phi \ll Y^{-\frac{1}{2} -J}.$$
Hence, by (\ref{k int near 0 w Z}), the integral in  (\ref{k int near 0 w error}) is now
\begin{equation}\label{k int near 0 w H_J}
\frac{\exp\left(\frac{k+1}{k^2}\zeta(\frac{k+1}{k})\Gamma(\frac{1}{k})X^{\frac{1}{k}} - \frac{1}{2}\right)}{(2\pi)^{\frac{k+2}{2}}X^{\frac{3}{2}}Y^{\frac{1}{2}}}  \left( \int _0^Z \mathfrak{R} \left(\exp\left(-\phi+H_J(\phi) \right)\right)\phi^{-1/2} \, d\phi + O(Y^{-J}) \right).
\end{equation}

We have 
$$\exp(H_J(\phi)) = \sum_{j=0}^\infty \frac{H_J(\phi)^j}{j!}.$$
The method of estimation of $\mathfrak{R}(H(\phi))$ can be used again to show that $|H_J(\phi)| \ll Y^{-\frac{1}{2}}\phi^{\frac{3}{2}}$.  Thus $|H_J(\phi)| \le Y^{-\frac{1}{4}}$.   
Therefore
$$\left| \int_0^Z e^{-\phi} \mathfrak{R}(H_J(\phi)^j)\phi^{-\frac{1}{2}}\, d\phi \right| \le Y^{-\frac{j}{4}}$$
and so 
$$\sum_{j=4J+4}^\infty \int_0^Z e^{-\phi} \mathfrak{R}\left(\frac{H_J(\phi)^j}{j!}\right)\phi^{-\frac{1}{2}}\, d\phi \ll Y^{-J-1}.$$

We are left to deal with 
\begin{equation} \label{truncated 0 to Z}
\int_0^Z e^{-\phi}\sum_{j=0}^{4J+3} \mathfrak{R}\left(\frac{H_J(\phi)^j}{j!}\right)\phi^{-\frac{1}{2}}\, d\phi.
\end{equation}
If $Y$ is fixed, then $H_J(\phi)$ is a polynomial in $i\phi^{1/2}$ of degree $2J+2$ with real coefficients.   Moreover, the coefficient of $(i\phi^{1/2})^j$ in $H_J(\phi)$ is given by 
$$(a_j + b_j Y^{-1})  Y^{1-\frac{j}{2}} = a_j (Y^{-\frac{1}{2}})^{j-2}+ b_j (Y^{-\frac{1}{2}})^{j},$$
which is itself a real polynomial in $Y^{-\frac{1}{2}}$ of degree $j$ with a zero of order $j-2$.  
The expression $\sum_{j=0}^{4J+3} \left(\frac{H_J(\phi)^j}{j!}\right)$ is therefore a real polynomial in $i\phi^\frac{1}{2}$ of degree at most $L = (2J+2)(4J+3)$.  This polynomial can be written as
$$\sum_{h=0}^L p_h(Y^{-\frac{1}{2}})(i\phi^\frac{1}{2})^h$$
where the coefficients $p_h(z)$ are polynomials in $z$ of degree at most $h$. In particular, we have $p_0(z) = 1$,  ${p_1(z) \equiv p_2(z) \equiv 0}$, and for $h\ge3$,  $p_h(0)=0$.  The polynomials $p_h(z)$ have parity that agrees with the parity of $h$.

For $0\le h\le L$, we have
$$\int_Z^\infty e^{-\phi}\phi^\frac{h-1}{2}\, d\phi \le Y^{-J} \int_0^\infty e^{-\phi}\phi^\frac{h-1}{2}\, d\phi \ll Y^{-J}.$$
Therefore, by (\ref{k int near 0 w H_J}),
\begin{align} \nonumber
&\int_{-3/(8\pi X)}^{3/(8\pi X)} \rho^{-n} \exp(\Phi_k (\rho e(\Theta))-2\pi i n \Theta)\, d\Theta 
\\& = \frac{\exp\left(\frac{k+1}{k^2}\zeta(\frac{k+1}{k})\Gamma(\frac{1}{k})X^{\frac{1}{k}} - \frac{1}{2}\right)}{(2\pi)^{\frac{k+2}{2}}X^{\frac{3}{2}}Y^{\frac{1}{2}}} \left( I + O(Y^{-J})\right)
\end{align}
where
\begin{align}
I & = \int _0^Z \mathfrak{R} \left(\exp\left(-\phi+H_J(\phi) \right)\right)\phi^{-\frac{1}{2}} \, d\phi \nonumber \\ &
= \sum_{\substack{h=0\\h\text{ even}}}^L  p_h(Y^{-\frac{1}{2}}) \int _0^Z e^{-\phi}\phi^{\frac{h-1}{2}} \, d\phi + O(Y^{-J}) \nonumber \\ &
= \Gamma(\frac{1}{2}) + \sum_{h=2}^{L/2} \Gamma\left(h+\frac{1}{2}\right)  p_{2h}(Y^{-\frac{1}{2}}) + O(Y^{-J}).\label{c_j polynomial}
\end{align}
Recall that $p_{2h}$ is an even polynomial, so the sum above is in fact a polynomial in $Y^{-1}$.  Let $c_j$ denote the coefficient of $Y^{-j}$ in that polynomial.  Then
\begin{align}
& \int_{-3/(8\pi X)}^{3/(8\pi X)} \rho^{-n} \exp(\Phi_k (\rho e(\Theta))-2\pi i n \Theta)\, d\Theta = \nonumber \\ &
\frac{\exp\left(\frac{k+1}{k^2}\zeta(\frac{k+1}{k})\Gamma(\frac{1}{k})X^{\frac{1}{k}} - \frac{1}{2}\right)}{(2\pi)^{\frac{k+2}{2}}X^{\frac{3}{2}}Y^{\frac{1}{2}}} \left( \pi^{\frac{1}{2}} + \sum_{j=1}^{J-1} c_jY^{-j} + O(Y^{-J})\right).
\label{main term}\end{align}

If we replace $J$ by $J+1$, this is the expression given in (\ref{main result}).  The remainder of the proof consists of showing that 
\begin{align*} &\int_{\mathcal{U} \setminus[-3/(8\pi X),3/(8\pi X)]} \rho^{-n} \exp(\Phi_k (\rho e(\Theta))-2\pi i n \Theta)\, d\Theta
\\& \ll \frac{\exp\left(\frac{k+1}{k^2}\zeta(\frac{k+1}{k})\Gamma(\frac{1}{k})X^{\frac{1}{k}} - \frac{1}{2}\right)}{(2\pi)^{\frac{k+2}{2}}X^{\frac{3}{2}}Y^{\frac{1}{2}}} Y^{-J}.
\end{align*}
Suppose that $\Theta\in\mathfrak{M}(1,0) \setminus [-3/(8\pi X), 3/(8\pi X)]$.  Then $\Delta$, given by (\ref{Delta}), is less than $\frac{4}{5}$.  Applying Lemma \ref{major arcs lemma} with $q=1, a=0$, we see that 
\begin{align*}
\mathfrak R \Phi_k(\rho e(\Theta)) & 
\le \Gamma(\frac{k+1}{k})\zeta(\frac{k+1}{k}) (X\Delta)^\frac{1}{k} + O(X^{\frac{1}{2}}|\theta|^{\frac{1}{2}}\log X) \\ & 
 \le \frac{1}{k} \Gamma(\frac{1}{k})\zeta(\frac{k+1}{k}) \left(\frac{4X}{5}\right)^\frac{1}{k} + O(X^{\frac{1}{2k}+\varepsilon}).
 \end{align*}
 Note that for $\Theta\in\mathfrak{M}(1,0)$ we have $|\theta| = |\Theta| \le X^{1/k-1}$.
Therefore, 
\begin{align}\label{rest of M(1,0)}
&\left|\int_{\mathfrak{M}(1,0) \setminus [-3/(8\pi X), 3/(8\pi X)]} \rho^{-n} \exp(\Phi_k (\rho e(\Theta))-2\pi i n \Theta)\, d\Theta\right| 
\\ &\ll \frac{\exp\left(\frac{k+1}{k^2}\zeta(\frac{k+1}{k})\Gamma(\frac{1}{k})X^{\frac{1}{k}} - \frac{1}{2}\right)}{(2\pi)^{\frac{k+2}{2}}X^{\frac{3}{2}}Y^{\frac{1}{2}}}Y^{-J}.\nonumber
\end{align}

We next study the integral on the remaining major arcs.  Suppose that $\Theta\in\mathfrak{M}(q,a)$ with $q>1$.  Then we have $q\le X^{1/k}$ and $\theta = \Theta - a/q$ satisfies $|\theta|\le q^{-1}X^{1/k-1}$.  Thus 
$$q^{1/2+\varepsilon}\log X\left(1+ X^\frac{1}{2}|\theta|^\frac{1}{2}\right) \ll X^{\frac{1}{2k}+\varepsilon},$$
and by Lemma \ref{major arcs lemma} we have 
$$\Phi_k(e^{-1/X}e(\Theta)) = \Gamma\left(\frac{k+1}{k}\right)\left(\frac{X}{1-2\pi i X\Theta}\right)^{1/k}\sum_{j=1}^\infty \frac{S_k(q_j, a_j)}{j^{\frac{k+1}{k}}q_j} + O\left(X^{1/(2k)+\varepsilon}\right)$$
In order to obtain an estimate for the integral over the major arcs, we will first need a bound on the sum $$\sum_{j=1}^\infty\frac{|S_k(q_j, a_j)|}{j^{\frac{k+1}{k}}q_j}.$$  

If $q\mid j$, then we have $q = (q,j)$, i.e. $q_j = 1$ and $S_k(q_j, a_j) = 1$.   
 On the other hand, if $q\nmid j$,  then $q_j>1$ and Lemma \ref{exponential sum lemma} tells us that there is a constant $\delta_k>0$ such that $|S_k(q_j, a_j)| \le q_j(1-\delta_k)$.  Thus the sum satisfies
 \begin{align*}
 \sum_{j=1}^\infty\frac{|S_k(q_j, a_j|)}{j^{\frac{k+1}{k}}q_j} &= \sum_{\substack{j=1\\ q\nmid j}}^\infty\frac{|S_k(q_j, a_j)|}{j^{\frac{k+1}{k}}q_j} + \sum_{\substack{j=1\\ q\mid j}}^\infty\frac{|S_k(q_j, a_j)|}{j^{\frac{k+1}{k}}q_j}
 \le \sum_{\substack{j=1\\ q\nmid j}}^\infty \frac{1-\delta_k}{j^{\frac{k+1}{k}}} + \sum_{\substack{j=1\\ q\mid j}}^\infty\frac{1}{j^{\frac{k+1}{k}}}
\\ &
 = (1-\delta_k)\left(1-\frac{1}{q^{\frac{k+1}{k}}}\right) \zeta\left(\frac{k+1}{k}\right) + \frac{1}{q^{\frac{k+1}{k}}}\zeta\left(\frac{k+1}{k}\right) \\ &
 = \left(1-\delta_k + \frac{\delta_k}{q^{\frac{k+1}{k}}} \right) \zeta(\frac{k+1}{k}) 
 < \left(1- \frac{\delta_k}{2}\right) \zeta(\frac{k+1}{k}).
 \end{align*}
 Here we used the readily verifiable fact that $\sum_{q\mid j} j^{-\alpha} = q^{-\alpha}\zeta(\alpha)$.
 We now have
 $$\left|\Re \left(\Gamma(\frac{k+1}{k})\left(\frac{X}{1-2\pi i X\Theta}\right)^{1/k}\sum_{j=1}^\infty \frac{S_k(q_j, a_j)}{j^{\frac{k+1}{k}}q_j} \right)\right|\le \frac{1}{k}(1- \frac{\delta_k}{2})\zeta(\frac{k+1}{k})\Gamma(\frac{1}{k}) X^{1/k}  .$$
 Let $\widetilde{\mathfrak{M}} = \mathfrak{M}\setminus\mathfrak{M}(1,0).$
 Then the above argument proves that there is a constant $\delta$ with $0<\delta<1$ such that
 \begin{align} 
 \int_{\widetilde{\mathfrak{M}}} \rho^{-n} \exp(\Phi_k (\rho e(\Theta))-2\pi i n \Theta)\, d\Theta & 
 \ll \exp\left(\frac{n}{X} + \frac{\delta}{k}\zeta(\frac{k+1}{k})\Gamma(\frac{1}{k}) X^{1/k} \right) \nonumber\\ & 
 \ll \frac{\exp\left(\frac{k+1}{k^2}\zeta(\frac{k+1}{k})\Gamma(\frac{1}{k})X^{\frac{1}{k}} - \frac{1}{2}\right)}{(2\pi)^{\frac{k+2}{2}}X^{\frac{3}{2}}Y^{\frac{1}{2}}} Y^{-J}.\label{major arc estimate}
 \end{align}

Finally, we deal with the minor arcs.  We will need one more estimate for $\Phi_k(\rho e(\Theta))$.

\begin{lemma}\label{minor arcs lemma} Let $\Theta\in\mathfrak{m}$.  Then, with all definitions as above,
$$\Phi_k(\rho e(\Theta)) \ll X^{\frac{1}{k} + \varepsilon - \frac{1}{k2^{k-1}}} .$$
\end{lemma}
\noindent{\bf Proof}~ Let $K\in\mathbb{Z}$ be a parameter at our disposal.  As in the proof of Lemma \ref{major arcs lemma}, we have
\begin{equation}\label{minor arcs first step}
\Phi_k(e^{-1/X}e(\Theta)) = \sum_{j=1}^K \frac{1}{j}  \int_1^\infty jk x^{k-1} X^{-1} e^{-x^kj/X} \sum_{n\le x} e(jn^k\Theta)\,dx + O\left(\left(\frac{X}{K}\right)^{1/k}\right).
\end{equation}
For each $j$, we use Dirichlet's Theorem to choose $a_j\in\mathbb{Z}, q_j\in\mathbb{N}$, so that 
$$\left| j\Theta - \frac{a_j}{q_j}\right| \le q_j^{-1}X^{\frac{1}{k}-1},\quad q_j \le X^{1-\frac{1}{k}}.$$
We now use Weyl's inequality \cite{Vaughan1997} to obtain
$$\sum_{n\le x} e(jn^k\Theta) \ll x^{1+\varepsilon - 2^{-(k-1)}} + x^{1+\varepsilon}q_j^{-2^{-(k-1)}} + x^{1+\varepsilon}\left(\frac{q_j}{x^k}\right)^{2^{-(k-1)}}.$$
Note that for any $\lambda>0$,
$$\int_1^\infty x^\lambda\left( jk x^{k-1} X^{-1} e^{-x^kj/X}\right)\,dx \ll\left( \frac{X}{j}\right)^{\lambda/k}.
$$
So, the main term of (\ref{minor arcs first step}) is 
$$\ll  \sum_{j=1}^K \frac{1}{j}  \left(\left( \frac{X}{j}\right)^{\frac{1+\varepsilon}{k} - \frac{1}{k2^{k-1}}} +\left( \frac{X}{j}\right)^{\frac{1+\varepsilon}{k}}q_j^{-2^{-(k-1)}} + \left( \frac{X}{j}\right)^{\frac{1+\varepsilon}{k} - \frac{1}{2^{k-1}}} q_j^{2^{-(k-1)}}\right).$$

Since $q_j \le X^{1-\frac{1}{k}}$, the last term is 
$$ \left( \frac{X}{j}\right)^{\frac{1+\varepsilon}{k} - \frac{1}{2^{k-1}}} q_j^{2^{-(k-1)}} \le  \left( \frac{X}{j}\right)^{\frac{1+\varepsilon}{k} }  \left( \frac{X}{j^k}\right)^{- \frac{1}{k2^{k-1}}} .$$
As $\Theta\not\in\mathfrak{M}$, we must have $jq_j > X^{1/k}$.  So we can replace the middle term as well:
$$\left( \frac{X}{j}\right)^{\frac{1+\varepsilon}{k}}q_j^{-2^{-(k-1)}} < \left( \frac{X}{j}\right)^{\frac{1+\varepsilon}{k}}\left( \frac{X}{j^k}\right)^{- \frac{1}{k2^{k-1}}}.$$
Putting all of this together, we see that
\begin{align*}& \Phi_k(e^{-1/X}e(\Theta))  \\ &
\ll X^{\frac{1+\varepsilon}{k} - \frac{1}{k2^{k-1}}} \sum_{j=1}^K\left( \left( \frac{1}{j}\right)^{1+\frac{1+\varepsilon}{k} - \frac{1}{k2^{k-1}}} + \left( \frac{1}{j}\right)^{1+\frac{1+\varepsilon}{k} - \frac{1}{2^{k-1}}} \right) + \left(\frac{X}{K}\right)^{1/k} \\ &
\ll X^{\frac{1}{k} + \varepsilon - \frac{1}{k2^{k-1}}}  + \left(\frac{X}{K}\right)^{1/k}. \\ &
\end{align*}
Letting $K\rightarrow\infty$ gives the desired result.\hfill $\square$\newline\newline
For $\Theta\in\mathfrak{m}$, we have by Lemma \ref{minor arcs lemma}
$$\Phi_k(\rho e(\Theta)) \ll X^{\frac{1}{k}-\frac{1}{k^2}}.$$ 
Therefore
\begin{equation}\label{minor arc estimate}
\int_{\mathfrak{m}} \rho^{-n} \exp(\Phi_k (\rho e(\Theta))-2\pi i n \Theta)\, d\Theta \ll \frac{\exp\left(\frac{k+1}{k^2}\zeta(\frac{k+1}{k})\Gamma(\frac{1}{k})X^{\frac{1}{k}} - \frac{1}{2}\right)}{(2\pi)^{\frac{k+2}{2}}X^{\frac{3}{2}}Y^{\frac{1}{2}}} Y^{-J}.
\end{equation}

Combining (\ref{rest of M(1,0)}), (\ref{major arc estimate}), and (\ref{minor arc estimate}) completes the proof of Theorem \ref{main theorem}.

\section{Proof of Theorem \ref{difference theorem}} 
We will use the notation from the previous section.  Recall that 
$$p^k(n) = \int_{\mathcal U} \rho^{-n} \exp(\Phi_k (\rho e(\Theta))-2\pi i n \Theta)\, d\Theta.$$
Hence we have
$$p^k(n+1) - p^k(n) = \int_{\mathcal U} \rho^{-n} \exp(\Phi_k (\rho e(\Theta))-2\pi i n \Theta)(\rho^{-1}e^{-2\pi i \Theta} -1)\, d\Theta.$$
Since $|\rho^{-1}e^{-2\pi i \Theta} -1| \le e^{1/X}+1\le 4$, the contribution from $|\Theta|>3/(8\pi X)$ is 
$$\ll  \frac{\exp\left(\frac{k+1}{k^2}\zeta(\frac{k+1}{k})\Gamma(\frac{1}{k})X^{\frac{1}{k}} - \frac{1}{2}\right)}{(2\pi)^{\frac{k+2}{2}}X^{\frac{3}{2}}Y^{\frac{1}{2}}} X^{-2}Y^{-J}$$
by the proof of Theorem \ref{main theorem}.  On the other hand, when $|\Theta|\le3/(8\pi X)$, we have
$$\rho^{-1}e^{-2\pi i \Theta} -1 = \exp\left(\frac{1}{X} - 2\pi i \Theta\right) - 1 = \frac{1}{X} -2\pi i \Theta + O(X^{-2}).$$
We thus deduce that
\begin{align}\label{difference first step}
&p^k(n+1) - p^k(n) = -2\pi i \int_{-3/(8\pi X)}^{3/(8\pi X)} \rho^{-n} \exp(\Phi_k (\rho e(\Theta))-2\pi i n \Theta) \Theta\, d\Theta \\ &
+  \frac{\exp\left(\frac{k+1}{k^2}\zeta(\frac{k+1}{k})\Gamma(\frac{1}{k})X^{\frac{1}{k}} - \frac{1}{2}\right)}{(2\pi)^{\frac{k+2}{2}}X^{\frac{3}{2}}Y^{\frac{1}{2}}} \left( \pi^{\frac{1}{2}} + \sum_{j=1}^{J-1} c_jY^{-j} + O(Y^{-J})\right)\left(X^{-1} + O(X^{-2})\right). \nonumber
\end{align}

It remains to evaluate 
$$-2\pi i \int_{-3/(8\pi X)}^{3/(8\pi X)} \rho^{-n} \exp(\Phi_k (\rho e(\Theta))-2\pi i n \Theta) \Theta\, d\Theta.$$
The methods here are similar to those of Theorem \ref{main theorem}, so we only outline the major differences here.  The extra factor of $\Theta$ in the integrand means that (\ref{k int near 0 w error}) becomes 
\begin{align}\label{diff int near 0 w error}&
 \int_{-3/(8\pi X)}^{3/(8\pi X)}\rho^{-n} \exp(\Xi_k(\rho e(\Theta))-2\pi i n \Theta)\Theta\, d\Theta 
 \\ & + O\left(X^{-3/2} \exp\left(\frac{k+1}{k^2}\zeta(\frac{k+1}{k})\Gamma(\frac{1}{k})X^{\frac{1}{k}} - \frac{1}{2} -\delta X^{\frac{1}{k}}\right)\right) \nonumber
 \end{align}
and the righthand side of (\ref{k int near 0 w G}) becomes 
\begin{equation}
\frac{\exp\left(\frac{k+1}{k^2}\zeta(\frac{k+1}{k})\Gamma(\frac{1}{k})X^{\frac{1}{k}} - \frac{1}{2}\right)}{(2\pi)^{\frac{k}{2}}X^{\frac{1}{2}}}\int_{-3/(8\pi X)}^{3/(8\pi X)} \exp\left(-Y(2\pi X \Theta)^2 + G(\Theta) \right)\Theta\, d\Theta, \label{diff int near 0 w G}
\end{equation}
where $G(\Theta)$ is defined by (\ref{G defn}).

We rewrite the integral in (\ref{diff int near 0 w G}) as
\begin{align*}
&\int_{0}^{3/(8\pi X)}\left( \exp( G(\Theta)) -  \exp(G(-\Theta))\right)\exp\left(-Y(2\pi X \Theta)^2 \right)\Theta\, d\Theta \\ &
= 2i \Im \int_{0}^{3/(8\pi X)} \exp\left(G(\Theta) -Y(2\pi X \Theta)^2 \right)\Theta \, d\Theta \\ &
= \frac{i}{4\pi^2 X^2 Y}  \Im \int_{0}^{9Y/16} \exp\left(H(\phi) -\phi \right) \, d\phi,
\end{align*}
using the change of variables $\phi = (2\pi X\Theta)^2 Y$ with $H$ defined by (\ref{H defn}).
Putting this into (\ref{diff int near 0 w error}), we see that 
\begin{align}
& -2\pi i \int_{-3/(8\pi X)}^{3/(8\pi X)} \rho^{-n} \exp(\Phi_k (\rho e(\Theta))-2\pi i n \Theta) \Theta\, d\Theta \\ &
= \frac{\exp\left(\frac{k+1}{k^2}\zeta(\frac{k+1}{k})\Gamma(\frac{1}{k})X^{\frac{1}{k}} - \frac{1}{2}\right)}{(2\pi)^{\frac{k}{2}+1}X^{\frac{5}{2}}Y}\left(\Im\int_{0}^{9Y/16} \exp\left(H(\phi) -\phi \right) \, d\phi + O(Y^{-J})\right).
\end{align}
The computations from this point on are entirely similar to Theorem \ref{main theorem}, and we end up with \begin{align} & 
\Im\int_{0}^{9Y/16} \exp\left(H(\phi) -\phi \right) \, d\phi  \nonumber \\ &
= \int_0^Z e^{-\phi}\sum_{j=0}^{4J+3} \Im\left(\frac{H_J(\phi)^j}{j!}\right)\, d\phi + O(Y^{-J}) \nonumber \\ &
= \sum_{\substack{h=1\\h\text{ odd}}}^L  p_h(Y^{-\frac{1}{2}}) \int _0^Z e^{-\phi}\phi^{\frac{h}{2}} \, d\phi + O(Y^{-J}) \nonumber \\ &
= \sum_{h=2}^{L/2}  \Gamma\left(h+\frac{1}{2}\right)p_{2h-1}(Y^{-\frac{1}{2}})  + O(Y^{-J})\label{odd polynomial}
\end{align}

Recall that $p_{2h-1}$ is an odd polynomial in $Y^{-1/2}$.  For $j\ge 1$, let $\widetilde{c}_j$ denote the coefficient of $(Y^{-1/2})^{2j-1}$ in the polynomial in (\ref{odd polynomial}).  Then we have
\begin{align}
& -2\pi i \int_{-3/(8\pi X)}^{3/(8\pi X)} \rho^{-n} \exp(\Phi_k (\rho e(\Theta))-2\pi i n \Theta) \Theta\, d\Theta \nonumber \\ &
= \frac{\exp\left(\frac{k+1}{k^2}\zeta(\frac{k+1}{k})\Gamma(\frac{1}{k})X^{\frac{1}{k}} - \frac{1}{2}\right)}{(2\pi)^{\frac{k}{2}+1}X^{\frac{5}{2}}Y^{\frac{1}{2}}}\left( \sum_{j=1}^{J}  \widetilde{c}_jY^{-j} + O(Y^{-J-1})\right). \nonumber
\end{align}
Combining this with (\ref{difference first step}) and letting $d_j = c_j+\widetilde{c}_j$, we have
\begin{equation}
p^k(n+1) - p^k(n) =   \frac{\exp\left(\frac{k+1}{k^2}\zeta(\frac{k+1}{k})\Gamma(\frac{1}{k})X^{\frac{1}{k}} - \frac{1}{2}\right)}{(2\pi)^{\frac{k+2}{2}}X^{\frac{5}{2}}Y^{\frac{1}{2}}} \left( \pi^{\frac{1}{2}} + \sum_{j=1}^{k-1} d_jY^{-j} + O(Y^{-k})\right). \nonumber
\end{equation}
The error term of $Y^{-k}$ is due to the error of $X^{-1}$ in (\ref{difference first step}).

\section{Computing the Coefficients} \label{computing the coefficients}

In certain applications, it may be useful to know the values of the coefficients $c_j$, which appear in (\ref{main result}).   As seen in (\ref{c_j polynomial}), these coefficients satisfy
\begin{equation}\label{coefficients}
\sum_{j=1}^{L} c_j Y^{-j} = \sum_{h=2}^{L/2} \Gamma\left(h+ \frac{1}{2}\right) p_{2h}(Y^{-1/2}),
\end{equation} 
where the polynomials $p_{2h}$ are given by
\begin{equation}\label{p polynomials}
\sum_{h=0}^{L} p_{2h}(Y^{-1/2}) \phi^{h} = \Re\sum_{\ell=0}^{4J+3} \frac{H_J(\phi)^\ell}{\ell!}
\end{equation}
with 
\begin{equation} \label{H definition}
H_J(\phi) = \sum_{m=3}^{2J+2} i^m (a_m + b_m Y^{-1}) \phi^{m/2} Y^{1-m/2}
\end{equation}
where
$$ a_m = \left(\frac{2k^2}{k+1} \right)\frac{\Gamma(m+\frac{1}{k}) }{m!\ \Gamma(\frac{1}{k})}, \quad b_m = \frac{1}{2}\left(\frac{2k^2}{k+1} \right)\frac{\Gamma(m+\frac{1}{k}) }{m!\ \Gamma(\frac{1}{k})} - \frac{1}{2m}.$$  

Given a particular value of $J$, one could input (\ref{H definition}) and (\ref{p polynomials}) into a computer algebra program, and expand it out into powers of $\phi$ to obtain an expression for each of the polynomials $p_h(z)$, which in turn could be put into (\ref{coefficients}) to obtain $c_j$.  Note that $L = (4J+3)(2J+2)$ is the degree of $\sum_{\ell=0}^{4J+3} \frac{H_J(\phi)^\ell}{\ell!}$ as a polynomial in $\phi^{1/2}$.   Because of the error term of $Y^{-J-1}$ in (\ref{main result}), it is only useful to find $c_j$ up to $j = J$.  However, it is necessary to compute all of the polynomials $p_h$, $h=0,\ldots, L$, in order to compute any of the coefficients $c_j$.  For $J = 1$, there are $29$ polynomials to compute to obtain the first coefficient $c_1 = -\frac{\sqrt{\pi}}{24k^2}(k^2+\frac{5}{2}k+1)$.

\subsection*{Acknowledgements} 
This work was completed at Pennsylvania State University as part of the author's Ph.D. dissertation research.  The author would like to thank her supervisor, R.C. Vaughan, for his guidance throughout the project.  Special thanks also to James Fan for helping with calculations in Maple. 

\bibliographystyle{amsplain}
\bibliography{library}

\end{document}